\documentclass[11pt]{article}
\usepackage{enumerate}
\usepackage{amssymb,a4wide,latexsym,makeidx,epsfig,fleqn}
\usepackage{amsthm}
\usepackage{amsmath}
\usepackage{enumerate}
\usepackage{indentfirst}
\newtheorem{theorem}{Theorem}[section]

\newtheorem{lemma}[theorem]{Lemma}

\newtheorem{corollary}[theorem]{Corollary}

\begin{document}
\textwidth 150mm \textheight 225mm
\title{$D^Q$-integral and $D^L$-integral generalized wheel graphs
\thanks{ Supported by the National Natural Science Foundation of China (No. 12271439).}}
\author{Yirui Chai$^{a,b}$, Ligong Wang$^{a,b}$\thanks{Corresponding author.}, Yuwei Zhou$^{a,b}$\\
	\author{\textbf{Ligong Wang \thanks{Corresponding author: lgwangmath@163.com}}}\\
	{\small $^a$ School of Mathematics and Statistics, }\\ {\small  Northwestern Polytechnical University, Xi'an, Shaanxi 710129, P.R. China. }\\
	{\small $^b$ Xi'an-Budapest Joint Research Center for Combinatorics,}\\{\small Northwestern Polytechnical University, Xi'an, Shaanxi 710129, P.R. China. }\\
	{\small E-mail: yiruichai@163.com, lgwangmath@163.com, yuweizhoumath@163.com} }
\date{}

\maketitle
\begin{center}
\begin{minipage}{120mm}
\vskip 0.3cm
\begin{center}
{\small {\bf Abstract}}
\end{center}
{\small A graph $ G $ is said to be $ M $-integral (resp. $A$-integral, $D$-integral, $D^L$-integral or $D^Q$-integral) if all eigenvalues of its matrix $ M $ (resp. adjacency matrix $ A(G) $, distance matrix $ D(G) $, distance Laplacian matrix $D^L(G)$ or distance signless Laplacian matrix $D^Q(G)$)  are integers. Lu et al. [Discrete Math, 346 (2023)] defined the generalized wheel graph $ GW(a,m,n) $ as the graph $aK_{m}\nabla C_{n}$, and obtained all $ D $-integral generalized wheel graphs $aK_{m} \nabla C_{n}.$ Based on the above research, in this paper, we determine all $D^L$-integral and $D^Q$-integral generalized wheel graphs $aK_{m} \nabla C_{n} $ respectively. As byproducts, we give a sufficient and necessary condition for the join of regular graphs $ G_{1} \nabla G_{2} $ to be $ D^{L} $-integral, from which we can get infinitely many new classes of $ D^{L} $-integral graphs according to the large number of research results about the $A$-integral graphs.

 }
\end{minipage}
\end{center}
\section{Introduction }
\label{sec:ch6-introduction}
Let $ G = (V(G),E(G)) $ be a simple, undirected and connected graph on $ n $ vertices and $ A(G) $ be the adjacency matrix of $ G, $ where $ V(G) $ is the vertex set and $ E(G) $ is the edge set. For a graph $ G $, the distance between two vertices $ u,v \in V(G) $, denoted by $ d_{G}(u,v) $ or $ d(u,v) $, is defined to be the length of the shortest path between $ u $ and $ v $. We denote by $ D(G)=(d(u,v))_{u,v \in V(G) } $ the distance matrix of $ G $, and by $ Tr(G) $, the transmission matrix of $ G $, the diagonal matrix of the row sums of $ D(G) $. In 2013, Aouchiche and Hansen \cite{AoHa} introduced the Laplacian and the signless Laplacian for the distance matrix of a connected graph, defined as $ D^{L}(G) = Tr(G)-D(G) $ and $ D^{Q}(G) = Tr(G)+D(G) $ respectively. As usual, if $ M $ is a real symmetric matrix associated to the graph $ G, $ then the graph $ G $ is called $ M $-integral when all eigenvalues of $ M $ are integers. The $ M $-spectrum of the graph $ G $ consists of all eigenvalues of its matrix $ M $ together with multiplicities.

A graph $ G $ is regular if every vertex has the same degree. The union $ G_{1} \cup G_{2} $ of two graphs $ G_{1} $ and $ G_{2} $ is the graph whose vertex set is $ V(G_{1}) \cup V(G_{2}) $ and the edge set is $ E(G_{1}) \cup E(G_{2}). $ We write $ kG $ for the union of $ k $ copies of $ G. $ The join $ G_{1} \nabla G_{2} $ of two disjoint graphs $ G_{1} $ and $ G_{2} $ is the graph obtained from $ G_{1} \cup G_{2} $ by adding all possible edges from the vertices of $ G_{1} $ to those in $ G_{2}. $ As usual, we denote the complete graph and cycle on $ n $ vertices by $ K_n $  and $ C_n. $ In this paper, we also use some basic concepts of number theory. We say $d$ divides $n$ and write $d \mid n$ whenever $n=c d$ for some $c$. If $ d $ divides two integers $ a $ and $ b $, then $ d $ is called a common divisor of $ a $ and $ b $. The greatest common divisor $ d $ of $ a $ and $ b $ and is denoted by $ d=(a, b) $. If $ (a, b) = 1 $ then $ a $ and $ b $ are said to be relatively prime. Unless otherwise stated, we use the standard notations and terminologies in \cite{Ap, BoMu, BrHa}.

In spectral graph theory, one of the most important questions is to characterize graphs for which all eigenvalues of a matrix associated to the graph are integers. In 1974, Harary and Schwenk \cite{HaSc} posed the notion of integral graphs, which sparked considerable interest and research on integral graphs among scholars. In 2002, Balali{\'{n}}ska et al. \cite{BaCvRaSiSt} demonstrated a survey of results on integral graphs. Over the past four decades, the exploration of integral graphs has been a significant focus in research. We refer the interested reader to the surveys \cite{AhBeMo, BuCv, ChFeYuZh, LuHuHu, MuKl, Roi, Wa, WaLiHo, WaLiZh} and so on for more results.

In regard to the $ D $-integral graph, in 2010, Ili{\'{c}} \cite{Il} characterized the distance spectra of integral circulant graphs and proved that these graphs are $ D $-integral. In 2011, Renteln \cite{Re} discovered that the absolute order graphs of the Coxeter groups are $ D $-integral. In 2015, Pokorn{\'{y}} et al. \cite{PoHiStMi} characterized $ D $-integral graphs in the classes of complete split graphs, multiple complete split-like graphs, extended complete split-like graphs and multiple extended complete split-like graphs, and showed that no nontrivial tree can be $ D $-integral. In 2015, Yang and Wang \cite{YaWa} gave a sufficient and necessary condition for the complete $ r $-partite graph $ {K_{{p_1},{p_2}, \cdots ,{p_r}}} \cong {K_{{a_1} \cdot {p_1},{a_2} \cdot {p_2}, \cdots ,{a_s} \cdot {p_s}}} $ to be $D$-integral and constructed infinitely many new classes of $ D $-integral graphs with $ s = 1,2,3,4. $ In 2016, H{\'{i}}c et al. \cite{HiPoTr} found infinitely many new classes of $ D $-integral complete $ r $-partite graph $ {K_{{p_1},{p_2}, \cdots ,{p_r}}} \cong {K_{{a_1} \cdot {p_1},{a_2} \cdot {p_2}, \cdots ,{a_s} \cdot {p_s}}} $ when $ s = 5,6 $. In 2021, Huang and Li \cite{HuLi1} gave some sufficient and necessary conditions for the Cayley graphs to be $D$-integral over generalized dihedral groups. In 2021, Huang and Li \cite{HuLi2} displayed the sufficient and necessary conditions for the Cayley graphs to be $D$-integral over abelian groups and dicyclic groups. In 2023, Mirafzal \cite{Mi} proved that the line graph of the crown graph $ L(Cr(n)) $ is $ D $-integral. In 2024, Wu et al. \cite{WuFeYu} presented some criteria for the distance integrality of quasiabelian 2-Cayley graphs.

With regard to $D^L$-integral and $D^Q$-integral graph, in 2016, Zhao et al. \cite{ZhLiMe} gave a sufficient and necessary condition for the complete $ r $-partite graph $ {K_{{p_1},{p_2}, \cdots ,{p_r}}} \cong {K_{{a_1} \cdot {p_1},{a_2} \cdot {p_2}, \cdots ,{a_s} \cdot {p_s}}} $ to be $D^Q$-integral and constructed infinitely many new classes of $ D^Q $-integral graphs with $ s = 1,2,3 $. In 2017, Da Silva Junior et al. \cite{DaDeDe} considered the $D^L$-integrality and $D^Q$-integrality of complete split graphs, multiple complete split-like graphs, extended complete split-like graphs and multiple extended complete split-like graphs based on Pokorn{\'{y}} \cite{PoHiStMi}. For more results, one can refer to \cite{BaRaSa, KaDe, Pa} and the review \cite{LiShXZh}.

Our inspiration for this paper comes from \cite{LuLiuLu}. In 2023, Lu et al. \cite{LuLiuLu} obtained all $ D $-integral generalized wheel graphs $aK_{m} \nabla C_{n}.$  Analogously, we determine all $D^L$-integral and $D^Q$-integral generalized wheel graphs $aK_{m} \nabla C_{n} $ respectively.
The rest of this paper is organized as follows.
In Section \ref{sec:ch-$D^{Q}$-integral-2}, we derive all $D^{Q}$-integral generalized wheel graphs $aK_{m} \nabla C_{n} $ which consist of an infinitely class of graphs and $ 17 $ scattered graphs.
In Section \ref{sec:ch-$D^{L}$-integral}, we give a useful sufficient and necessary condition for the join of two regular graphs $ G_{1} \nabla G_{2} $ to be $ D^{L} $-integral, from which we can get infinitely many new classes of $ D^{L} $-integral graphs according to the large number of research results of predecessors about the $A$-integral graphs. Furthermore, we determine all $D^{L}$-integral generalized wheel graphs which consist of three infinite class of graphs.


\section{$D^{Q}$-integral generalized wheel graphs $aK_{m} \nabla C_{n} $}
\label{sec:ch-$D^{Q}$-integral-2}

In this section, we shall provide the distance signless Laplacian spectrum of the generalized wheel graphs $ GW(a,m,n) $ and completely determine all $D^{Q}$-integral generalized wheel graphs.

In 2023, Lu et al. \cite{LuLiuLu} defined the generalized wheel graph $ GW(a,m,n) $ as the graph $aK_{m}\nabla C_{n}$, where $ K_m $ is the complete graph on $ m $ vertices and $ C_n $ is the cycle graph on $ n $ vertices. In 2017, Da Silva Junior et al. \cite{DaDeDe} determined the $D^{Q}$-characteristic polynomials for graphs $ G_{1} \nabla G_{2} $ where $G_{i}$ is $r_{i}$-regular, for $i = 1, 2$. Thus, we naturally obtain the following Lemma \ref{le:ch3-1} which is vital for what follows.
\noindent\begin{lemma}\label{le:ch3-1} (\cite{DaDeDe}) For $i = 1, 2$, let $G_{i}$ be an $r_{i}$-regular graph with $n_{i}$ vertices. If the eigenvalues of the adjacency matrix of $G_{i}$ are given by $ r_{i}=\lambda_{1}^{(i)}\geq \lambda_{2}^{(i)} \geq \cdots \geq \lambda_{n_{i}}^{(i)} $, then the distance signless Laplacian spectrum of $ G_{1} \nabla G_{2} $ consists of eigenvalues $ 2(n_{1}-2)+n_{2}-r_{1}-\lambda_{j}^{(1)} $ for $ j=2,3,\cdots,n_{1}, $ and $ 2(n_{2}-2)+n_{1}-r_{2}-\lambda_{j}^{(2)} $ for $ j=2,3,\cdots,n_{2}, $ the remaining two eigenvalues are
	$$ \frac{-8+5(n_{1}+n_{2})-2(r_{1}+r_{2})}{2}\pm \frac{\sqrt{(3(n_{1}-n_{2})-2(r_{1}-r_{2}))^2+4 n_{1}n_{2}}}{2}.$$
\end{lemma}

\noindent\begin{lemma}\label{le:ch4-2} (\cite{BrHa}) The adjacency spectrum of $ aK_{m} $ consists of the eigenvalues $ m-1 $ and $ -1 $ with multiplicities $ a $ and $ a(m-1) $ respectively, the adjacency spectrum of $ C_{n} $ is $ \{ 2\cos(\frac{2\pi j}{n});0\leq j \leq n-1 \} $.
\end{lemma}

Using Lemmas \ref{le:ch3-1} and \ref{le:ch4-2}, we can get the distance signless Laplacian spectrum of the generalized wheel graph $ GW(a,m,n) $ on $ am+n $ vertices as follow.

\noindent\begin{theorem}\label{th:ch4-1} The distance signless Laplacian spectrum of the generalized wheel graph $ GW(a,m,n) $ consists of the eigenvalues $ 2(a-1)m+n-2 $, $ (2a-1)m+n+2 $ and $ am+2n-6-2\cos(\frac{2\pi j}{n}) $ with multiplicities $ a-1 $, $ a(m-1) $ and $ 1 $ respectively, where $ 1 \leq j \leq n-1, $ the remaining two eigenvalues are $$ \frac{(5a-2)m+5n-10}{2}\pm\frac{\sqrt{((3a-2)m-3n+6)^2+4amn}}{2}. $$
\end{theorem}

By Theorem \ref{th:ch4-1}, the following result is immediate.

\noindent\begin{corollary}\label{cr:ch4-1} The distance signless Laplacian spectrum of $ GW(1,m,n) $ on $ m+n $ vertices, consists of the eigenvalues $ m+n-2 $ with multiplicities $ m-1 $ and $ m+2n-6-2\cos(\frac{2\pi j}{n}), $ where $ 1 \leq j \leq n-1, $ the remaining two eigenvalues are $$ \frac{3m+5n-10}{2}\pm\frac{\sqrt{(m-3n+6)^2+4mn}}{2}. $$
\end{corollary}
Next, we shall completely determine all $D^{Q}$-integral the generalized wheel graphs $ GW(a,m,n) $. Now we start by considering the case where $ a=1, $ i.e., the generalized wheel graph $ GW(1,m,n)=K_{m}\nabla C_{n} $. Lemmas \ref{le:ch-3} is useful for the following proofs.

\noindent\begin{lemma}\label{le:ch-3} For $ \forall $ $ x \in \mathbb{N} $, if $ x \equiv y $ (mod 2), then  $ x \equiv y^2 $ (mod 2).	
\end{lemma}
\noindent {\bf Proof.} Since $ x^2-x=x(x-1) \equiv 0 $ (mod 2) for $ \forall $ $ x \in \mathbb{N} $, we know that $ x \equiv x^2 $ (mod 2).	By $ x \equiv x^2 $ (mod 2) and $ x^2 \equiv y^2 $ (mod 2), we have the conclusion $ x \equiv y^2 $ (mod 2). \hfill$\qedsymbol$
%
%

\noindent\begin{lemma}\label{le:ch4-4} The generalized wheel graph $ GW(1,m,n) $ on $ m+n $ vertices is $D^{Q}$-integral if and only if one of the following cases holds
	\begin{enumerate}[(i)]
		\item $ n=3,$  $m \geq 1 $.
		
		\item $ n=4,$  $m = 5 $.
		
		\item $ n=6,$  $m = 5$.
		
		\item $ n=6,$  $m = 9$.
		
		\item $ n=6,$  $m = 16$.	
		
		\item $ n=6,$  $m = 35$.
	\end{enumerate}	
\end{lemma}
\noindent {\bf Proof.} By Corollary \ref{cr:ch4-1}, we have the distance signless Laplacian spectrum of $ GW(1,m,n) $ consists of the eigenvalues $ m+n-2 $ with multiplicities $ m-1 $ and $ m+2n-6-2\cos(\frac{2\pi j}{n}), $ where $ 1 \leq j \leq n-1, $ the remaining two eigenvalues are
\begin{equation}\label{eq:ch4-1}
\frac{3m+5n-10}{2}\pm\frac{\sqrt{(m-3n+6)^2+4mn}}{2}.
\end{equation}

It is clear that $ 3m+5n-10 \equiv m-3n+6 $ (mod 2). Using Lemma \ref{le:ch-3} we have $ 3m+5n-10 \equiv (m-3n+6)^2+4mn $ (mod 2). Hence, (\ref{eq:ch4-1}) is integers if and only if $ (m-3n+6)^2+4mn $ is a perfect square. Moreover, it is apparent that $ \cos(\frac{2\pi j}{n}) $  is integral for any $ 1 \leq j \leq n-1 $ if and only if $ n \in \{3,4,6\}. $  Therefore, the generalized wheel graph $ GW(1,m,n) $ on $ m+n $ vertices is $D^{Q}$-integral if and only if $ (m-3n+6)^2+4mn $ is a perfect square and $ n \in \{3,4,6\} $.

The proof of sufficiency is straightforward by basic calculations. Thus $ GW(1,m,n) $ is obviously $D^{Q}$-integral when condition $ (i) $, $ (ii) $, $ (iii) $, $ (iv) $, $ (v) $ or $ (vi) $ holds. Next, we will consider the necessity.
Let $ t=(m-3n+6)^2+4mn=c^2 $, we discuss the following three cases.

{\bf Case 1.} $ n=3. $

In this case, we have $ t=(m-3)^2+12m=(m+3)^2. $ Therefore, $ t $ is a perfect square for any $ m \geq 1, $ that is, $GW(1,m,3)=K_{m}\nabla C_{3}$ is always $D^{Q}$-integral for any $ m \geq 1. $ $ (i) $ holds.

{\bf Case 2.} $ n=4. $

In this case, we have $ t=(m-6)^2+16m=(m+2)^2+32=c^2, $ that is, $ 32=[c+(m+2)][c-(m+2)]. $ Since 32 = 32 $\times$ 1 = 16 $\times$ 2 = 8 $\times$ 4. For the case $[c+(m+2)][c-(m+2)]= 32 \times 1$, we have $ 2c=33 $ which contradicts $ c \in \mathbb{Z} $. For the other two cases we have $ c=9 $ and $ m+2=7, $ or $ c=6 $ and $ m+2=2. $ Therefore, we obtain $ (ii) $ $ c=9 $ and $ m=5 $ from the former, and the other contradicts $ m \geq 1 $ and $ m \in \mathbb{Z}. $

{\bf Case 3.} $ n=6. $

In this case, we have $ t=(m-12)^2+24m=m^2+144=c^2, $ that is, $ 144=(c+m)(c-m). $ It follows that $ (iii) $ $ c=13 $ and $ m=5, $ $ (iv) $ $ c=15 $ and  $ m=9, $ $ (v) $ $ c=20 $ and $ m=16, $ $ (vi) $ $ c=37 $ and $ m=35. $ \hfill$\qedsymbol$

\hspace*{\fill}

The next result characterizes a sufficient and necessary condition for the generalized wheel graph $ GW(a,m,n) $ to be $D^{Q}$-integral with $ a \geq 2 $ and $ n \geq 3 $.

\noindent\begin{lemma}\label{le:ch4-5} For any integers $ a,m,n $ with $ a \geq 2 $ and $ n \geq 3 $, the generalized wheel graph $ GW(a,m,n)$ is $D^{Q}$-integral if and only if the positive integers $ a,m,n $ satisfy one of the following cases
	\begin{enumerate}[(i)]
		\item $ n=3 $ and $$ m=\frac{\alpha^2+6\alpha(a-2)-72a(a-1)}{2\alpha(3a-2)^2} $$
		for some $ \alpha \in \mathbb{N} $ with $ \alpha \geq \sqrt{72a(a-1)}, $ or
		
		$$ m=\frac{-\alpha^2+6\alpha(a-2)+72a(a-1)}{2\alpha(3a-2)^2} $$
		for some $ \alpha \in \mathbb{N} $ with $ \sqrt{72a(a-1)} \leq \alpha \leq 12(a-1). $
		
		\item $ n=4 $ and $$ m=\frac{\alpha^2+4\alpha(5a-6)-32a(7a-6)}{2\alpha(3a-2)^2} $$
		for some $ \alpha \in \mathbb{N} $ with $ \alpha \geq \sqrt{32a(7a-6)}, $ or
		
		$$ m=\frac{-\alpha^2+4\alpha(5a-6)+32a(7a-6)}{2\alpha(3a-2)^2} $$
		for some $ \alpha \in \mathbb{N} $ with $ \sqrt{32a(7a-6)} \leq \alpha < 4(7a-6). $
		
		\item $ n=6 $ and $$ m=\frac{\alpha^2+48\alpha(a-1)-144a(5a-4)}{2\alpha(3a-2)^2} $$
		for some $ \alpha \in \mathbb{N} $ with $ \alpha \geq \sqrt{144a(5a-4)}, $ or
		
		$$ m=\frac{-\alpha^2+48\alpha(a-1)+144a(5a-4)}{2\alpha(3a-2)^2} $$
		for some $ \alpha \in \mathbb{N} $ with $ \sqrt{144a(5a-4)} \leq \alpha < 12(5a-4). $
	\end{enumerate}	
\end{lemma}
\noindent {\bf Proof.}  By Theorem \ref{th:ch4-1}, we have the distance signless Laplacian spectrum of of the generalized wheel graph $ GW(a,m,n) $ consists of the eigenvalues $ 2(a-1)m+n-2 $, $ (2a-1)m+n+2 $ and $ am+2n-6-2\cos(\frac{2\pi j}{n}) $ with multiplicities $ a-1 $, $ a(m-1) $ and $ 1 $, respectively, where $ 1 \leq j \leq n-1, $ the remaining two eigenvalues are
\begin{equation}\label{eq:ch4-2}
\frac{(5a-2)m+5n-10}{2}\pm\frac{\sqrt{((3a-2)m-3n+6)^2+4amn}}{2}.
\end{equation}

It is clear that $ (5a-2)m+5n-10 \equiv (3a-2)m-3n+6 $ (mod 2). Using Lemma \ref{le:ch-3} we have $ (5a-2)m+5n-10 \equiv ((3a-2)m-3n+6)^2+4amn $ (mod 2). Hence, (\ref{eq:ch4-2}) is integers if and only if $ ((3a-2)m-3n+6)^2+4amn $ is a perfect square. Moreover, it is apparent that $ \cos(\frac{2\pi j}{n}) $  is integral for any $ 1 \leq j \leq n-1 $ if and only if $ n \in \{3,4,6\} $. Therefore, the generalized wheel graph $ GW(a,m,n) $ is $D^{Q}$-integral if and only if $ ((3a-2)m-3n+6)^2+4amn $ is a perfect square and $ n \in \{3,4,6\} $.

The proof of sufficiency is straightforward by basic calculations. Thus the generalized wheel graph $ GW(a,m,n) $ is obviously $D^{Q}$-integral when condition $ (i) $, $ (ii) $, or $ (iii) $ holds. Next, we will consider the necessity.
Let $ t=(3a-2)m-3n+6)^2+4amn $ and $ t=c^2 $, we discuss the following three cases.

{\bf Case 1.} $ n=3. $

In this case, we have $ t=((3a-2)m-3)^2+12am=(3a-2)^2 m^2-6(a-2)m+9>9 $ and $ c>3 $. Moreover, we have $ (3a-2)^2 m^2-6(a-2)m+9=c^2, $ that is, $ (3a-2)^2 m^2-6(a-2)m+9-c^2=0 $. It follows that $ m=\frac{3(a-2)\pm\sqrt{(3a-2)^2c^2-72a(a-1)}}{(3a-2)^2}. $  Let
\begin{equation}\label{eq:ch4-3}
m=\frac{3(a-2) \pm p}{(3a-2)^2} ,
\end{equation} where
\begin{equation}\label{eq:ch4-4}
(3a-2)^2c^2-72a(a-1)=p^2,
\end{equation}
for some $ p\geq 0 $ such that $ m,p \in \mathbb{Z} $. Under this circumstance, we observe that
\begin{align} \label{eq:ch4-5}
72a(a-1)&=(3a-2)^2c^2-p^2  \nonumber    \\
~&=[(3a-2)c+p][(3a-2)c-p] \nonumber    \\
~&=\alpha[(3a-2)c-p] \nonumber    \\
~&=(3a-2)\alpha c-\alpha p,
\end{align}
where
\begin{equation}\label{eq:ch4-6}
\alpha=(3a-2)c+p.
\end{equation}
Combining Equations (\ref{eq:ch4-5}) and (\ref{eq:ch4-6}), we obtain $ c=\frac{\alpha^2+72a(a-1)}{2(3a-2)\alpha} $. By Equation (\ref{eq:ch4-4}), we have $ p=\frac{\alpha^2-72a(a-1)}{2\alpha} $. Since $ p \geq 0 $, we have $ \alpha \geq \sqrt{72a(a-1)} $. Combining Equation (\ref{eq:ch4-3}), it is easy to verify that  $ m=\frac{\alpha^2+6\alpha(a-2)-72a(a-1)}{2\alpha(3a-2)^2} $ or $ m=\frac{-\alpha^2+6\alpha(a-2)+72a(a-1)}{2\alpha(3a-2)^2}. $

Since $ m > 0 $, when $ m=\frac{\alpha^2+6\alpha(a-2)-72a(a-1)}{2\alpha(3a-2)^2}, $ we have $ \alpha > 6a $, it is easy to verify that $ 6a<\sqrt{72a(a-1)} $ when $ a \geq 2 $. In this instance, we have  $ \alpha \geq \sqrt{72a(a-1)} $.

On the other hand, when $ m=\frac{-\alpha^2+6\alpha(a-2)+72a(a-1)}{2\alpha(3a-2)^2}, $ we have $ \alpha < 12(a-1) $ according to $ m > 0 $. In this instance, we have  $ \sqrt{72a(a-1)} \leq \alpha < 12(a-1) $.
We complete the proof of $ (i) $.

{\bf Case 2.} $ n=4. $

In this case, we have $ t=((3a-2)m-6)^2+16am=(3a-2)^2 m^2-4(5a-6)m+36 $. Moreover, we have $ (3a-2)^2 m^2-4(5a-6)m+36=c^2, $ that is, $ (3a-2)^2 m^2-4(5a-6)m+36-c^2=0 $. It follows that $ m=\frac{2(5a-6)\pm\sqrt{(3a-2)^2c^2-32a(7a-6)}}{(3a-2)^2}. $ Let
\begin{equation}\label{eq:ch4-7}
m=\frac{2(5a-6) \pm p}{(3a-2)^2} ,
\end{equation} where
\begin{equation}\label{eq:ch4-8}
(3a-2)^2c^2-32a(7a-6)=p^2,
\end{equation}
for some $ p\geq 0 $ such that $ m,p \in \mathbb{Z} $. Under this circumstance, we observe that
\begin{align} \label{eq:ch4-9}
32a(7a-6)&=(3a-2)^2c^2-p^2  \nonumber    \\
~&=[(3a-2)c+p][(3a-2)c-p] \nonumber    \\
~&=\alpha[(3a-2)c-p] \nonumber    \\
~&=(3a-2)\alpha c-\alpha p,
\end{align}
where
\begin{equation}\label{eq:ch4-10}
\alpha=(3a-2)c+p.
\end{equation}
Combining Equations (\ref{eq:ch4-9}) and (\ref{eq:ch4-10}), we obtain $ c=\frac{\alpha^2+32a(7a-6)}{2(3a-2)\alpha} $. By Equation (\ref{eq:ch4-8}), we have $ p=\frac{\alpha^2-32a(7a-6)}{2\alpha} $. Since $ p \geq 0 $, we have $ \alpha \geq \sqrt{32a(7a-6)} $. Combining Equation (\ref{eq:ch4-7}), it is easy to verify that  $  m=\frac{\alpha^2+4\alpha(5a-6)-32a(7a-6)}{2\alpha(3a-2)^2} $ or $ m=\frac{-\alpha^2+4\alpha(5a-6)+32a(7a-6)}{2\alpha(3a-2)^2} $.

Since $ m > 0 $, when $ m=\frac{\alpha^2+4\alpha(5a-6)-32a(7a-6)}{2\alpha(3a-2)^2} $, we have $ \alpha > 8a $, it is easy to verify that $ 8a<\sqrt{32a(7a-6)} $ when $ a \geq 2 $. In this instance, we have  $ \alpha \geq \sqrt{32a(7a-6)} $.

On the other hand, when $ m=\frac{-\alpha^2+4\alpha(5a-6)+32a(7a-6)}{2\alpha(3a-2)^2} $, we have $ \alpha < 4(7a-6) $ according to $ m > 0 $. In this instance, we have  $ \sqrt{32a(7a-6)} \leq \alpha < 4(7a-6) $.
We complete the proof of $ (ii) $.

{\bf Case 3.} $ n=6. $

In this case, we have $
t=((3a-2)m-12)^2+24am=(3a-2)^2 m^2-48(a-1)m+144. $ Moreover, we have $ (3a-2)^2 m^2-48(a-1)m+144=c^2, $ that is, $ (3a-2)^2 m^2-48(a-1)m+144-c^2=0. $ It follows that $ m=\frac{24(a-1)\pm\sqrt{(3a-2)^2c^2-144a(5a-4)}}{(3a-2)^2}. $ Let
\begin{equation}\label{eq:ch4-11}
m=\frac{24(a-1)\pm p}{(3a-2)^2} ,
\end{equation} where
\begin{equation}\label{eq:ch4-12}
(3a-2)^2c^2-144a(5a-4)=p^2,
\end{equation}
for some $ p\geq 0 $ such that $ m,p \in \mathbb{Z} $. Under this circumstance, we observe that
\begin{align} \label{eq:ch4-13}
144a(5a-4)&=(3a-2)^2c^2-p^2  \nonumber    \\
~&=[(3a-2)c+p][(3a-2)c-p] \nonumber    \\
~&=\alpha[(3a-2)c-p] \nonumber    \\
~&=(3a-2)\alpha c-\alpha p,
\end{align}
where
\begin{equation}\label{eq:ch4-14}
\alpha=(3a-2)c+p.
\end{equation}
Combining Equations (\ref{eq:ch4-13}) and (\ref{eq:ch4-14}), we obtain $ c=\frac{\alpha^2+144a(5a-4)}{2(3a-2)\alpha} $. By Equation (\ref{eq:ch4-12}), we have $ p=\frac{\alpha^2-144a(5a-4)}{2\alpha} $. Since $ p \geq 0 $, we have $ \alpha \geq \sqrt{144a(5a-4)} $. Combining Equation (\ref{eq:ch4-11}), it is easy to verify that  $   m=\frac{\alpha^2+48\alpha(a-1)-144a(5a-4)}{2\alpha(3a-2)^2} $ or $ m=\frac{-\alpha^2+48\alpha(a-1)+144a(5a-4)}{2\alpha(3a-2)^2} $.

Since $ m > 0 $, when $ m=\frac{\alpha^2+48\alpha(a-1)-144a(5a-4)}{2\alpha(3a-2)^2} $, we have $ \alpha > 12a $, it is easy to verify that $ 12a <\sqrt{144a(5a-4)} $ when $ a \geq 2 $. In this instance, we have  $ \alpha \geq \sqrt{144a(5a-4)} $.

On the other hand, when $ m=\frac{-\alpha^2+48\alpha(a-1)+144a(5a-4)}{2\alpha(3a-2)^2} $, we have $ \alpha < 12(5a-4) $
according to $ m > 0 $. In this instance, we have  $ \sqrt{144a(5a-4)} \leq \alpha < 12(5a-4) $.
We complete the proof of $ (iii) $.\hfill$\qedsymbol$

\noindent\begin{lemma}\label{le:ch-6} (\cite{Ap} (Euclid's lemma)) If $ a \mid bc $ and $ (a, b)=1 $, then $ a \mid c $.
\end{lemma}
\noindent {\bf Proof.} Since $ (a, b)=1 $ we can write $ 1=ax+by $. Therefore $ c=acx+bcy $. But $ a \mid a c x $ and $ a \mid bcy $, so $ a \mid c $. \hfill$\qedsymbol$

\hspace*{\fill}

Using Lemmas \ref{le:ch4-5} and \ref{le:ch-6}, we further characterize a more specific condition for the generalized wheel graph $ GW(a,m,n) $ to be $D^{Q}$-integral with $ a \geq 2 $ and $ n \geq 3 $ as follow.

\noindent\begin{lemma}\label{le:ch4-7} For positive integers $ a \geq 2 $, $ m \geq 1 $, and and $ n \geq 3 $, if the generalized wheel graph $ GW(a,m,n) $ is $D^{Q}$-integral, one of the following conditions is satisfied:
	\begin{enumerate}[(i)]
		\item $ n=3,$  $m \leq 2 $.
		
		\item $ n=4,$  $m \leq 8 $.
		
		\item $ n=6,$  $m \leq 31$.
	\end{enumerate}	
\end{lemma}

\noindent {\bf Proof.} By Lemma \ref{le:ch4-5}, we will consider the following three cases.

{\bf Case 1.} $ n=3. $

We consider the following two subcases by Lemma \ref{le:ch4-5}.

{\bf Subcase 1.1.} $ m=\frac{\alpha^2+6\alpha(a-2)-72a(a-1)}{2\alpha(3a-2)^2} $ for some $ \alpha \in \mathbb{N} $ with $ \alpha \geq \sqrt{72a(a-1)}. $

In this case, it leads to $ \alpha \mid [\alpha^2+6\alpha(a-2)-72a(a-1)] $ and $ 2 \alpha \mid [\alpha^2+6\alpha(a-2)-72a(a-1)] $ according to $ m \in \mathbb{Z}. $ The former leads to $ \alpha \mid 72a(a-1) $ according to $ \alpha \mid [\alpha^2+6\alpha(a-2)]. $ If $ \alpha $ is odd then $ (\alpha,2)=1. $ But $ \alpha \mid 72a(a-1), $ so we have $ \alpha \mid 36a(a-1) $ by Lemma \ref{le:ch-6}. If $ \alpha $ is even then $ 2\alpha \mid 72a(a-1) $ according to $ 2\alpha \mid [\alpha^2+6\alpha(a-2)], $ that is, $ \alpha \mid 36a(a-1). $ Therefore, the condition $ \alpha \mid 36a(a-1) $ always holds, we are easily to get $ \alpha \leq 36a(a-1) $ according to $ 36a(a-1) \geq 0. $

Suppose to the contrary that $ m \geq 3. $ Then $ \alpha^2+6\alpha(a-2)-72a(a-1) \geq 6 \alpha(3a-2)^2 $ and so $ \alpha \geq 27a^2 - 39a + 18 + 3(3a-2) \sqrt{9a^2-14a+9} $ or $ \alpha \leq 27a^2 - 39a + 18 - 3(3a-2) \sqrt{9a^2-14a+9}. $ From the former solution, it is easy to verify that
\begin{equation}
\begin{aligned}
\alpha &\geq 27a^2 - 39a + 18 + 3(3a-2) \sqrt{9a^2-14a+9}  \nonumber    \\
~&>27a^2 - 39a + 18 + 3(3a-2)(3a-\frac{7}{3}) \nonumber    \\
~&=54a^2-78a+32 \nonumber    \\
~&>36a(a-1) \quad (a \geq 2),
\end{aligned}
\end{equation}
which contradicts $ \alpha \leq 36a(a-1). $ Similarly, from the latter solution, it is easy to verify that
\begin{equation}
\begin{aligned}
\alpha &\leq 27a^2 - 39a + 18 - 3(3a-2) \sqrt{9a^2-14a+9}  \nonumber    \\
~&<27a^2 - 39a + 18 - 3(3a-2)(3a-\frac{7}{3}) \nonumber    \\
~&=4 < \sqrt{72a(a-1)} \quad (a \geq 2),
\end{aligned}
\end{equation}
which contradicts $ \alpha \geq \sqrt{72a(a-1)}. $

{\bf Subcase 1.2.} $ m=\frac{-\alpha^2+6\alpha(a-2)+72a(a-1)}{2\alpha(3a-2)^2} $ for some $ \alpha \in \mathbb{N} $ with $ \sqrt{72a(a-1)} \leq \alpha \leq 12(a-1). $

Suppose to the contrary that $ m \geq 3. $ Then $ -\alpha^2+6\alpha(a-2)+72a(a-1) \geq 6 \alpha(3a-2)^2 $ and so $ -27a^2 + 39a - 18 - 3(3a-2) \sqrt{9a^2-14a+9} \leq \alpha \leq -27a^2 + 39a - 18 + 3(3a-2) \sqrt{9a^2-14a+9}. $ From the right-hand side, it is easy to verify that
\begin{equation}
\begin{aligned}
\alpha &\leq -27a^2 + 39a - 18 + 3(3a-2) \sqrt{9a^2-14a+9} \nonumber    \\
~&<-27a^2 + 39a - 18 + 3(3a-2)(3a-\frac{5}{3}) \nonumber    \\
~&=6a-8 < \sqrt{72a(a-1)} \quad (a \geq 2),
\end{aligned}
\end{equation}
which contradicts $ \alpha \geq \sqrt{72a(a-1)}. $

Combining the above arguments, we have that $ m \leq 2 $ when $ n=3. $ Thus, $ (i) $ holds.

%
{\bf Case 2.} $ n=4. $

We consider the following two subcases by Lemma \ref{le:ch4-5}.

{\bf Subcase 2.1.} $ m=\frac{\alpha^2+4\alpha(5a-6)-32a(7a-6)}{2\alpha(3a-2)^2} $ for some $ \alpha \in \mathbb{N} $ with $ \alpha \geq \sqrt{32a(7a-6)}. $

In this case, it leads to $ \alpha \mid [\alpha^2+4\alpha(5a-6)-32a(7a-6)] $ and $ 2 \alpha \mid [\alpha^2+4\alpha(5a-6)-32a(7a-6)] $ according to $ m \in \mathbb{Z}. $ The former leads to $ \alpha \mid 32a(7a-6) $ according to $ \alpha \mid [\alpha^2+4\alpha(5a-6)]. $ If $ \alpha $ is odd then $ (\alpha,2)=1. $ But $ \alpha \mid 32a(7a-6), $ so we have $ \alpha \mid 16a(7a-6) $ by Lemma \ref{le:ch-6}. If $ \alpha $ is even then $ 2\alpha \mid 32a(7a-6) $ according to $ 2\alpha \mid [\alpha^2+4\alpha(5a-6)], $ that is, $ \alpha \mid 16a(7a-6). $ Therefore, the condition $ \alpha \mid 16a(7a-6) $ always holds, we are easily to get $ \alpha \leq 16a(7a-6) $ according to $ 16a(7a-6) \geq 0. $

Suppose to the contrary that $ m \geq 9. $ Then $ \alpha^2+4\alpha(5a-6)-32a(7a-6) \geq 18 \alpha(3a-2)^2 $ and so $ \alpha \geq 81a^2 - 118a + 48 + 3(3a-2) \sqrt{81a^2-128a+64} $ or $ \alpha \leq 81a^2 - 118a + 48 - 3(3a-2) \sqrt{81a^2-128a+64}. $ From the former solution, it is easy to verify that

\begin{equation}
\begin{aligned}
\alpha &\geq 81a^2 - 118a + 48 + 3(3a-2) \sqrt{81a^2-128a+64}  \nonumber    \\
~&>81a^2 - 118a + 48 + 3(3a-2)(9a- \frac{64}{9}) \nonumber    \\
~&=162a^2-236a+\frac{272}{3} \nonumber    \\
~&>16a(7a-6) \quad (a \geq 2),
\end{aligned}
\end{equation}
which contradicts $ \alpha \leq 16a(7a-6). $ Similarly, from the latter solution, it is easy to verify that
\begin{equation}
\begin{aligned}
\alpha &\leq 81a^2 - 118a + 48 - 3(3a-2) \sqrt{81a^2-128a+64}  \nonumber    \\
~&<81a^2 - 118a + 48 - 3(3a-2)(9a- \frac{64}{9}) \nonumber    \\
~&=\frac{16}{3} <\sqrt{32a(7a-6)} \quad (a \geq 2),
\end{aligned}
\end{equation}
which contradicts $ \alpha \geq \sqrt{32a(7a-6)}. $

{\bf Subcase 2.2.} $ m=\frac{-\alpha^2+4\alpha(5a-6)+32a(7a-6)}{2\alpha(3a-2)^2} $ for some $ \alpha \in \mathbb{N} $ with $ \sqrt{32a(7a-6)} \leq \alpha < 4(7a-6). $

Suppose to the contrary that $ m \geq 9. $ Then $ -\alpha^2+4\alpha(5a-6)+32a(7a-6) \geq 18 \alpha(3a-2)^2 $ and so $ -81a^2 + 118a - 48 - 3(3a-2) \sqrt{81a^2-128a+64} \leq \alpha \leq -81a^2 + 118a - 48 + 3(3a-2) \sqrt{81a^2-128a+64}. $ From the right-hand side, it is easy to verify that
\begin{equation}
\begin{aligned}
\alpha &\leq -81a^2 + 118a - 48 + 3(3a-2) \sqrt{81a^2-128a+64} \nonumber    \\
~&<-81a^2 + 118a - 48 + 3(3a-2)(9a- \frac{58}{9}) \nonumber    \\
~&=6a-\frac{28}{3} <\sqrt{32a(7a-6)} \quad (a \geq 2),
\end{aligned}
\end{equation}
which contradicts $ \alpha \geq \sqrt{32a(7a-6)}. $

Combining the above arguments, we have that $ m \leq 8 $ when $ n=4. $ Thus, $ (ii) $ holds.

{\bf Case 3.} $ n=6. $

We consider the following two subcases by Lemma \ref{le:ch4-5}.

{\bf Subcase 3.1.} $ m=\frac{\alpha^2+48\alpha(a-1)-144a(5a-4)}{2\alpha(3a-2)^2} $ for some $ \alpha \in \mathbb{N} $ with $ \alpha \geq \sqrt{144a(5a-4)}. $

In this case, it leads to $ \alpha \mid [\alpha^2+48\alpha(a-1)-144a(5a-4)] $ and $ 2 \alpha \mid [\alpha^2+48\alpha(a-1)-144a(5a-4)] $ according to $ m \in \mathbb{Z}. $ The former leads to $ \alpha \mid 144a(5a-4) $ according to $ \alpha \mid [\alpha^2+48\alpha(a-1)]. $ If $ \alpha $ is odd then $ (\alpha,2)=1. $ But $ \alpha \mid 144a(5a-4), $ so we have $ \alpha \mid 72a(5a-4) $ by Lemma \ref{le:ch-6}. If $ \alpha $ is even then $ 2\alpha \mid 144a(5a-4) $ according to $ 2\alpha \mid [\alpha^2+48\alpha(a-1)], $ that is, $ \alpha \mid 72a(5a-4). $ Therefore, the condition $ \alpha \mid 72a(5a-4) $ always holds, we are easily to get $ \alpha \leq 72a(5a-4) $ according to $ 72a(5a-4) \geq 0. $

Suppose to the contrary that $ m \geq 32. $ Then $ \alpha^2+48\alpha(a-1)-144a(5a-4) \geq 64 \alpha(3a-2)^2 $ and so $ \alpha \geq 4(72a^2 - 102a + 38 + (3a-2) \sqrt{576a^2-864a+361}) $ or $ \alpha \leq 4(72a^2 - 102a + 38 - (3a-2) \sqrt{576a^2-864a+361}). $ From the former solution, it is easy to verify that
\begin{equation}
\begin{aligned}
\alpha &\geq 4(72a^2 - 102a + 38 + (3a-2) \sqrt{576a^2-864a+361})  \nonumber    \\
~&>4(72a^2 - 102a + 38 + (3a-2)(24a-18)) \nonumber    \\
~&=4(144a^2-204a+74) \nonumber    \\
~&>72a(5a-4) \quad (a \geq 2),
\end{aligned}
\end{equation}
which contradicts $ \alpha \leq 72a(5a-4). $ Similarly, from the latter solution, it is easy to verify that
\begin{equation}
\begin{aligned}
\alpha &\leq 4(72a^2 - 102a + 38 - (3a-2) \sqrt{576a^2-864a+361})  \nonumber    \\
~&<4(72a^2 - 102a + 38 - (3a-2)(24a-18)) \nonumber    \\
~&=8 <\sqrt{144a(5a-4)} \quad (a \geq 2),
\end{aligned}
\end{equation}
which contradicts $ \alpha \geq \sqrt{144a(5a-4)}. $

{\bf Subcase 3.2.} $ m=\frac{-\alpha^2+48\alpha(a-1)+144a(5a-4)}{2\alpha(3a-2)^2} $ for some $ \alpha \in \mathbb{N} $ with $ \sqrt{144a(5a-4)} \leq \alpha < 12(5a-4). $

Suppose to the contrary that $ m \geq 32. $ Then $ -\alpha^2+48\alpha(a-1)+144a(5a-4) \geq 64 \alpha(3a-2)^2 $ and so $  4(-72a^2 +102a -38 -(3a-2) \sqrt{576a^2-864a+361})  \leq \alpha \leq 4(-72a^2 +102a -38 +(3a-2) \sqrt{576a^2-864a+361}). $ From the right-hand side, it is easy to verify that
\begin{equation}
\begin{aligned}
\alpha &\leq 4(-72a^2 +102a -38 +(3a-2) \sqrt{576a^2-864a+361}) \nonumber    \\
~&<4(-72a^2 +102a -38 +(3a-2)(24a-17)) \nonumber    \\
~&=4(3a-4)<\sqrt{144a(5a-4)} \quad (a \geq 2),
\end{aligned}
\end{equation}
which contradicts $ \alpha \geq \sqrt{144a(5a-4)}. $

Combining the above arguments, we have that $ m \leq 31 $ when $ n=6. $ Thus, $ (iii) $ holds.\hfill$\qedsymbol$

\hspace*{\fill}

In what follows, we try to find all positive integral solutions $ (a,m,n) $ satisfying Lemmas \ref{le:ch4-5} and \ref{le:ch4-7}.

\noindent\begin{lemma}\label{le:ch4-8} Let $ a \geq 2 $ and $ n \geq 3. $ Then the generalized wheel graph $ GW(a,m,n) $ is $D^{Q}$-integral if and only if the ordered triple $ (a,m,n) \in S $, where
	\begin{equation}
	\begin{aligned}
	S=&\{ (2,1,3),(2,1,4),(3,1,4),(4,2,4),(3,4,4),(4,1,6), \nonumber    \\
	&(5,1,6),(11,1,6),(4,2,6),(2,3,6),(5,3,6),(2,8,6) \}.
	\end{aligned}
	\end{equation}
\end{lemma}
\noindent {\bf Proof.} In accordance with Theorem \ref{th:ch4-1}, the proof of sufficiency is straightforward by basic calculations. Thus the generalized wheel graph $ GW(a,m,n) $ is $D^{Q}$-integral when $ (a,m,n) \in S $. Next, we will consider the necessity.

By the proof of Lemma \ref{le:ch4-5} we have the generalized wheel graph $ GW(a,m,n)=aK_{m}\nabla C_{n} $ is $D^{Q}$-integral if and only if $ ((3a-2)m-3n+6)^2+4amn $ is a perfect square and $ n \in \{3,4,6\} $. Suppose that the graph $aK_{m}\nabla C_{n}$ is $D^{Q}$-integral, we discuss the following three cases.

{\bf Case 1.} $ n=3. $

In this case, we have $ t=(3a-2)^2 m^2-6(a-2)m+9=c^2 $ by the proof of Lemma \ref{le:ch4-5}. But $ m \leq 2 $ according to Lemma \ref{le:ch4-7}, so when $ m = 1, $ we have $ t=(3(a-1))^2+16=c^2, $ that is, $ 16=[c+3(a-1)][c-3(a-1)]. $ It follows that $ c=5 $ and $ 3(a-1)=3, $ or $ c=4 $ and $ 3(a-1)=0. $ Therefore, we obtain $ c=5 $ and $ a=2 $ from the former, and the other contradict $ a \geq 2 $ and $ a \in \mathbb{Z}. $

When $ m = 2, $ we have $ t=(6a-5)^2+24=c^2, $ that is, $ 24=[c+(6a-5)][c-(6a-5)]. $ It follows that $ c=7 $ and $ 6a-5=5, $ or $ c=5 $ and $ 6a-5=1, $  which all contradict $ a \geq 2 $ and $ a \in \mathbb{Z}. $

{\bf Case 2.} $ n=4. $

In this case, we have $ t=(3a-2)^2 m^2-4(5a-6)m+36=c^2 $ by the proof of Lemma \ref{le:ch4-5}. But $ m \leq 8 $ according to Lemma \ref{le:ch4-7}, so when $ m = 1, $ we have $ t=(3a-\frac{16}{3})^2+\frac{320}{9}=c^2, $ that is, $ 320=[3c+(9a-16)][3c-(9a-16)]. $ It follows that $ 3c=18 $ and $ (9a-16)=2, $ $ 3c=21 $ and $ (9a-16)=11 $ $ 3c=24 $ and $ (9a-16)=16, $ $ 3c=42 $ and $ (9a-16)=38, $ or $ 3c=81 $ and $ (9a-16)=25. $ Therefore, we obtain respectively $ c=7 $ and $ a=3, $ $ c=6 $ and $ a=2, $ and the others contradict $ a \in \mathbb{Z}. $

When $ m = 2, $ we have $ t=4((3a-\frac{11}{3})^2+\frac{104}{9})=c^2, $ that is, $ 416=[3c+2(9a-11)][3c-2(9a-11)]. $ It follows that $ 3c=54 $ and $ 2(9a-11)=50, $ $ 3c=105 $ and $ 2(9a-11)=103, $ $ 3c=30 $ and $ 2(9a-11)=22, $ or $ 3c=21 $ and $ 2(9a-11)=5. $ Therefore, we obtain $ c=18 $ and $ a=4 $ from the former, and the others contradict $ a \in \mathbb{Z}. $

When $ m = 3, $ we have $ t=(9a-\frac{28}{3})^2+\frac{512}{9}=c^2, $ that is, $ 512=[3c+(27a-28)][3c-(27a-28)]. $ It follows that $ 3c=129 $ and $ 27a-28=127, $ $ 3c=66 $ and $ 27a-28=62, $ $ 3c=36 $ and $ 27a-28=28, $ or $ 3c=24 $ and $ 27a-28=8, $ which all contradict $ a \in \mathbb{Z} $.

When $ m = 4, $ we have $ t=4((6a-\frac{17}{3})^2+\frac{152}{9})=c^2, $ that is, $ 608=[3c+2(18a-17)][3c-2(18a-17)]. $ It follows that $ 3c=78 $ and $ 2(18a-17)=74, $ $ 3c=153 $ and $ 2(18a-17)=151, $ $ 3c=42$ and $ 2(18a-17)=34, $ or $ 3c=27 $ and $ 2(18a-17)=11. $ Therefore, we obtain $ c=26 $ and $ a=3 $ from the former, and the others contradict $ a \in \mathbb{Z}. $

Likewise, we find there is no solution that satisfy the conditions when $ 5 \leq m \leq 8 . $

{\bf Case 3.} $ n=6. $

In this case, we have $ t=(3a-2)^2 m^2-48(a-1)m+144=c^2 $ by the proof of Lemma \ref{le:ch4-5}. But $ m \leq 31 $ according to Lemma \ref{le:ch4-7}, so when $ m = 1, $ we have $ t=(3a-10)^2+96=c^2, $ that is, $ 96=[c+(3a-10)][c-(3a-10)]. $ It follows that $ c=25 $ and $ 3a-10=23, $ $ c=11 $ and $ 3a-10=5, $ $ c=10 $ and $ (3a-10)=2, $ or $ c=14 $ and $ (3a-10)=10. $ Therefore, we obtain respectively $ c=25 $ and $ a=11, $ $ c=11 $ and $ a=5, $ $ c=10 $ and $ a=4, $ and the other contradicts $ a \in \mathbb{Z}. $

When $ m = 2, $ we have $ t=4((3(a-2))^2+28)=c^2, $ that is, $ 112=[c+6(a-2)][c-6(a-2)]. $ It follows that $ c=16 $ and $ 6(a-2)=12, $ $ c=29 $ and $ 6(a-2)=27, $ or $ c=11 $ and $ 6(a-2)=3. $ Therefore, we obtain $ c=16 $ and $ a=4 $ from the former, and the others contradict $ a \in \mathbb{Z}. $

When $ m = 3, $ we have $ t=(9a-14)^2+128=c^2, $ that is, $ 128=[c+9a-14][c-(9a-14)]. $ It follows that $ c=33 $ and $ 9a-14=31, $ $ c=12 $ and $ 9a-14=4, $ or $ c=18 $ and $ 9a-14=14. $ Therefore, we obtain respectively $ c=33 $ and $ a=5 $, $ c=12 $ and $ a=2 $, and the other contradicts $ a \in \mathbb{Z}. $

When $ m = 8, $ we have $ t=(24(a-1))^2+208=c^2, $ that is, $ 208=[c+24(a-1)][c-24(a-1)]. $ It follows that $ c=28 $ and $ 24(a-1)=24, $ $ c=53 $ and $ 24(a-1)=51, $ or $ c=17 $ and $ 24(a-1)=9. $ Therefore, we obtain $ c=28 $ and $ a=2 $ from the former, and the others contradict $ a \in \mathbb{Z}. $

Likewise, we find there is no solution that satisfy the conditions when $ 4 \leq m \leq 7 $ or $ 9 \leq m \leq 31. $ \hfill$\qedsymbol$

\hspace*{\fill}

Combining Lemmas \ref{le:ch4-4} and Lemma \ref{le:ch4-8}, we completely find all $D^{Q}$-integral generalized wheel graph $ GW(a,m,n) $ as follow.

\noindent\begin{theorem}\label{th:ch4-2} The generalized wheel graph $ GW(a,m,n) $ is $D^{Q}$-integral if and only if one of the following cases holds
	\begin{enumerate}[(i)]	
		\item $ a=1,$  $n=3 $ and $ m \geq 1. $
		
		\item the ordered triple $ (a,m,n) \in S $, where
		\begin{equation}
		\begin{aligned}
		S=&\{ (1,5,4),(1,5,6),(1,9,6),(1,16,6),(1,35,6),(2,1,3),(2,1,4),(3,1,4), \nonumber    \\
		&(4,2,4),(3,4,4),(4,1,6),(5,1,6),(11,1,6),(4,2,6),(2,3,6),(5,3,6),(2,8,6) \}.
		\end{aligned}
		\end{equation}
	\end{enumerate}
\end{theorem}

\section{Results on $D^{L}$--integrality for the join of regular graphs}
\label{sec:ch-$D^{L}$-integral}

In this section, we determine the $D^{L}$-spectrum of $ G_{1} \nabla G_{2}, $ where $ G_{i} $ is $ r_{i} $-regular for $ i=1,2, $ and naturally get a sufficient and necessary condition for the join of two regular graphs of $ G_{1} \nabla G_{2} $ to be $D^{L}$-integral graphs.

\noindent\begin{lemma}\label{le:ch5-1} (\cite{BrHa}) Let $ M $ be a square matrix of order $ n $ that can be written in blocks as
	
	\[M = \left[ {\begin{array}{*{20}{c}}
		{{M_{1,1}}}&{{M_{1,2}}}& \cdots &{{M_{1,k}}}\\
		{{M_{2,1}}}&{{M_{2,2}}}& \cdots &{{M_{2,k}}}\\
		\vdots & \vdots & \ddots & \vdots \\
		{{M_{k,1}}}&{{M_{k,2}}}& \cdots &{{M_{k,k}}}
		\end{array}} \right],\]
	where $ M_{i,j}, $ $ 1 \leq i, j \leq  k, $ is the $ n_{i} \times m_{j} $ matrix such that its lines have constant sum equal to $ c_{ij}. $ Let $ \overline{M} = [c_{ij} ]_{k \times k}. $ Then, the eigenvalues of $ \overline{M} $ are also eigenvalues of $ M $.
\end{lemma}

In the following results, we determine the $D^{L}$-spectrum of the join of two regular graphs $ G_{1} \nabla G_{2}, $ where $ G_{i} $ is $ r_{i} $-regular for $ i=1,2. $

\noindent\begin{theorem}\label{th:ch5-1} For $i = 1, 2$, let $ G_{i} $ be an $ r_{i} $-regular graph with $ n_{i} $ vertices. If the eigenvalues of the adjacency matrix of $ G_{i} $ are given by $ r_{i}=\lambda_{1}^{(i)} \geq \lambda_{2}^{(i)} \geq \cdots \geq \lambda_{n_{i}}^{(i)} $, then the distance Laplacian spectrum of $ G_{1} \nabla G_{2} $ consists of eigenvalues $ 2n_{1}+n_{2}-r_{1}+\lambda_{j}^{(1)} $ for $ j=2,3,\cdots,n_{1}, $ and $ 2n_{2}+n_{1}-r_{2}+\lambda_{j}^{(2)} $ for $ j=2,3,\cdots,n_{2}, $ the remaining two eigenvalues are $ n_{1}+n_{2} $ and $ 0. $
\end{theorem}

\noindent {\bf Proof.} For $ i = 1,2 $, let $ G_{i} $ be a $ r_{i} $-regular graph on $ n_{i} $ vertices. Then, the distance matrix $ D(G_{1} \nabla G_{2}) $ can be written as
\[D(G_{1} \nabla G_{2}) = \left[ {\begin{array}{*{20}{c}}
	{2({\mathbb{J}_{{n_1} \times {n_1}}} - {\mathbb{I}_{{n_1} \times {n_1}}}) - A({G_1})}&{{\mathbb{J}_{{n_1} \times {n_2}}}}\\
	{{\mathbb{J}_{{n_2} \times {n_1}}}}&{2({\mathbb{J}_{{n_2} \times {n_2}}} - {\mathbb{I}_{{n_2} \times {n_2}}}) - A({G_2})}
	\end{array}} \right],\]
where $ \mathbb{J}_{{n_1} \times {n_1}} $ is the matrix of order $ n_1 $ with all elements equal to 1, $ \mathbb{I}_{{n_1} \times {n_1}} $ is the identity matrix of order $ n_1 $.
It is easy to verify that
\[Tr( {{G_1}\nabla {G_2}} ) = \left[ {\begin{array}{*{20}{c}}
	{(2\left( {{n_1} - 1} \right) + {n_2} - {r_1}){\mathbb{I}_{{n_1} \times {n_1}}}}&{{\mathbb{O}_{{n_1} \times {n_2}}}}\\
	{{\mathbb{O}_{{n_2} \times {n_1}}}}&{(2\left( {{n_2} - 1} \right) + {n_1} - {r_2}){\mathbb{I}_{{n_2} \times {n_2}}}}
	\end{array}} \right].\]
Thus, the matrix $ D^{L}(G_{1} \nabla G_{2})=Tr( {{G_1}\nabla {G_2}} ) - D(G_{1} \nabla G_{2}) $ has the form
\[\begin{array}{l}
{D^L}({G_1}\nabla {G_2})  \\
=\left[ {\begin{array}{*{20}{c}}
	{(2{n_1} + {n_2} - {r_1}){\mathbb{I}_{{n_1} \times {n_1}}} + A({G_1}) - 2{\mathbb{J}_{{n_1} \times {n_1}}}}&{{-\mathbb{J}_{{n_1} \times {n_2}}}}\\
	{{-\mathbb{J}_{{n_2} \times {n_1}}}}&{(2{n_2} + {n_1} - {r_2}){\mathbb{I}_{{n_2} \times {n_2}}} + A({G_2}) - 2{\mathbb{J}_{{n_2} \times {n_2}}}}
	\end{array}} \right]
\end{array}.\]

As a regular graph, $ G_{1} $ has the all-one vector $ \textbf{1}_{n_{1}} $ as an eigenvector corresponding to eigenvalue $ r_{1}, $ while all other eigenvectors are orthogonal to $ \textbf{1}_{n_1}. $ Let $ \lambda_{ v_{1}} $ be an arbitrary eigenvalue of the adjacency matrix of $ G_{1} $ with corresponding eigenvector $ v_{1} $, such that $ A({G_1})v_{1}=\lambda_{ v_{1}}v_{1} $ and $ v_{1} \bot \textbf{1}_{n_{1}}, $ where $ v_{1} \in \mathbb{R}^{n_{1}}. $

Then the vector $ u_{1}  =[{v_{1}}\textbf{0}]^{T} \in \mathbb{R}^{n_{1}+n_{2}} $ satisfies $ D^{L}(G_{1} \nabla G_{2})u_{1}  = (2{n_1} + {n_2} - {r_1} + \lambda_{ v_{1}}) u_{1}, $ that is, $ 2{n_1} + {n_2} - {r_1} + \lambda_{ v_{1}} $ is an eigenvalue of $ D^{L}(G_{1} \nabla G_{2}). $

Similarly, let $ \lambda_{ v_{2}} $ be an arbitrary eigenvalue of the adjacency matrix of $ G_{2} $ with corresponding eigenvector $ v_{2} $, such that $ A({G_2})v_{2}=\lambda_{ v_{2}}v_{2} $ and $ v_{2} \bot \textbf{1}_{n_{2}}, $ where $ v_{2} \in \mathbb{R}^{n_{2}}. $ Then the vector $ u_{2}  =[\textbf{0}{v_{2}}]^{T} \in \mathbb{R}^{n_{1}+n_{2}} $ satisfies $ D^{L}(G_{1} \nabla G_{2})u_{2}  = (2{n_2} + {n_1} - {r_2} + \lambda_{ v_{2}}) u_{2}, $ that is, $ 2{n_2} + {n_1} - {r_2} + \lambda_{ v_{2}} $ is an eigenvalue of $ D^{L}(G_{1} \nabla G_{2}). $

The matrix $ M_L $ can be written as
\[{M_L} = \left[ {\begin{array}{*{20}{c}}
	{{n_2}}&{{-n_2}}\\
	{{-n_1}}&{{n_1}}
	\end{array}} \right],\]
which has eigenvalues $ n_{1}+n_{2} $ and $ 0. $ By Lemma \ref{le:ch5-1}, the eigenvalues of the matrix $ M_L $ are also eigenvalues of $ D^{L}(G_{1} \nabla G_{2}). $\hfill$\qedsymbol$

\hspace*{\fill}

Next we naturally provide a sufficient and necessary condition for the join of two regular graphs $ G_{1} \nabla G_{2} $ to be $ D^{L} $-integral. The following corollary is immediate by Theorem \ref{th:ch5-1}.

\noindent\begin{corollary}\label{cr:ch5-1} The graph $ G_{1} \nabla G_{2} $ is $ D^{L} $-integral if and only if $ G_{i} $ is $ A $-integral, where $ G_{i} $ is $ r_{i} $-regular for $ i=1,2. $
\end{corollary}

Corollary \ref{cr:ch5-1} establishes the relationship between the distance Laplacian integrality and the adjacency integrality of the graph $ G_{1} \nabla G_{2} $ where $ G_{i} $ is $ r_{i} $-regular for $ i=1,2, $ from which we can get infinitely many new classes of $ D^{L} $-integral graphs according to the large number of research results of predecessors about the $A$-integral graphs.

Using Corollary \ref{le:ch4-2} and Lemma \ref{cr:ch5-1}, we get a sufficient and necessary condition for the generalized wheel graph $ GW(a,m,n) $ to be $ D^{L} $-integral. The following corollary is immediate.

\noindent\begin{corollary}\label{cr:ch5-2} The generalized wheel graph $ GW(a,m,n) $ is $ D^{L} $-integral if and only if  positive integers $ a,m,n $ satisfy one of the following cases
		\begin{enumerate}[(i)]
		\item $ a\geq1$, $m\geq1 $, $ n=3$.
		
		\item $ a\geq1$, $m\geq1 $, $ n=4$.
		
		\item $ a\geq1$, $m\geq1 $, $ n=6$.
	\end{enumerate}	
\end{corollary}

\section*{Declaration of competing interest}

The authors declare that they have no known competing financial interests or personal relationships that could have
appeared to influence the work reported in this paper.

\section*{Data availability}

No data was used for the research described in the article.

%

\end{document}